\documentclass[12pt,leqno]{article}
\usepackage{amsmath, amsthm, amsfonts, amssymb, color}
\usepackage{mathrsfs}
\theoremstyle{definition}

\newcommand{\scr}[1]{\mathscr #1}
\definecolor{wco}{rgb}{0.5,0.2,0.3}

\numberwithin{equation}{section} \theoremstyle{remark}

\newcommand{\ua}{\uparrow}

\title{{\bf Functional Inequalities for  Particle Systems on Polish
Spaces}}

\author{{\bf Michael R\"ockner$^{1}$ and Feng-Yu Wang$^{2}$\footnote{Supported in part by the DFG through the Forschergruppe
``Spectral Analysis, Asymptotic Distributions and Stochastic
Dynamics'', the BiBoS Research Centre,  NNSFC(10121101), and
RFDP(20040027009); Feng-Yu Wang: wangfy@bnu.edu.cn} }\\
\footnotesize{$^1$ Fakult\"at F\"ur Mathematik, Universit\"at Bielefeld, D-33501 Bielefeld, Germany}\\
\footnotesize{$^2$ School of Mathematical Sciences, Beijing Normal
University, Beijing 100875, China}}
\begin{document}
\maketitle
\begin{abstract}
Various Poincar\'e-Sobolev type inequalities are studied for a reaction-diffusion model of particle systems on
 Polish spaces. The systems we consider consist of finite  particles which are killed or produced
at certain rates, while  particles in the system move  on the
Polish space interacting with one another (i.e. diffusion). Thus,
the corresponding  Dirichlet form, which we call
reaction-diffusion Dirichlet form, consists of two parts: the
diffusion part  induced by certain Markov processes on the product
spaces $E^n (n\ge 1)$ which determine the motion of particles, and
the reaction part induced by a $Q$-process on $\mathbb Z_+$ and a
sequence of reference probability measures, where the $Q$-process
determines the variation of the number of particles and the
reference measures describe the locations of newly produced
particles.
 We prove that the validity of  Poincar\'e and weak
Poincar\'e inequalities are essentially due to the pure reaction part,
i.e. either of these inequalities holds if and only if it holds for
the pure reaction Dirichlet form, or equivalently,  for the
corresponding $Q$-process. But under a mild condition, stronger
inequalities rely on both parts:
 the reaction-diffusion Dirichlet form satisfies a super Poincar\'e inequality  (e.g. the
log-Sobolev inequality) if and only if so do both the
corresponding $Q$-process and the diffusion part.
 Explicit estimates of constants in the inequalities   are derived.
 Finally, some specific examples are presented to illustrate the main results.
\end{abstract}
 \vskip 2cm

\def\R{\mathbb R} \def\Z{\mathbb Z} \def\ff{\frac} \def\ss{\sqrt}
\def\N{\mathbb N}\def\gg{\gamma}
\def\dd{\delta} \def\DD{\Delta} \def\vv{\varepsilon} \def\rr{\rho}
\def\<{\langle} \def\>{\rangle} \def\GG{\Gamma}
\def\ll{\lambda} \def\LL{\Lambda} \def\nn{\nabla} \def\pp{\partial}
\def\d{\text{\rm{d}}} \def\bb{\beta} \def\aa{\alpha} \def\D{\scr D}
\def\E{\scr E} \def\si{\sigma} \def\ess{\text{\rm{ess}}}
\def\beg{\begin} \def\beq{\begin{equation}}  \def\F{\scr F}
\def\Ric{\text{\rm{Ric}}} \def\Hess{\text{\rm{Hess}}}\def\B{\scr B}
\def\e{\text{\rm{e}}} \def\ua{\underline a} \def\OO{\Omega} \def\b{\mathbf b}
\def\oo{\omega}     \def\tt{\tilde} \def\Ric{\text{\rm{Ric}}}
\def\cut{\text{\rm{cut}}} \def\P{\mathbb P} \def\ifn{I_n(f^{\bigotimes n})}
\def\fff{f(x_1)\dots f(x_n)} \def\ifm{I_m(g^{\bigotimes m})} \def\ee{\varepsilon}
\def\pm{\pi_{\mu}}   \def\p{\mathbf{p}}   \def\ml{\mathbf{L}}
 \def\C{\scr C}      \def\aaa{\mathbf{r}}     \def\r{r}
\def\gap{\text{\rm{gap}}} \def\prr{\pi_{\mu,\varrho}}  \def\r{\mathbf r}
\def\Z{\mathbb Z} \def\vrr{\varrho}

\section{Introduction}

In this paper we consider interacting particle systems in continnum, say in $\R^d$ or, more generally, in a Polish space $E$. At any given time
we have finitely (but arbitrarily) many particles interacting with one another (called the diffusion in $E$). In addition, the system kills or produces particles at certain rates (called the reaction). We refer e.g.  to \cite{Chenbook} where the corresponding discretized model (on the lattice $\Z^d$ instead of $E$) is analyzed in detail. The main aim of this paper is to derive functional inequalities for the Dirichlet form corresponding to these systems which, as
is well-known, gives information about their long-time behaviour.

Let us consider a system of finite  particles on $E$ such that the  number of  particles behaves as  a Markov chain
on $\Z_+$ generated by a regular $Q$-matrix
 $Q:=(q_{ij})_{i,j\ge 0}$. Assume  that the $Q$-process is reversible
w.r.t.  a probability measure $\vrr:=\{\vrr_i>0: i\ge 0\}$; that is,

\
\newline $(H_1)\ \  \vrr_iq_{ij}=\vrr_jq_{ji},\ \ i,j\ge
0.$

\ \newline Since we also consider the locations of particles, the state space
of the underlying Markov process for the particle system is the
following finite mutiple configuration space:

$$\GG_0:=\bigg\{\sum_{i=1}^n\dd_{x_i}:\  x_i\in E, n\ge 0 \bigg\},$$
where $\dd_x$ is the Dirac measure at $x$ and
$\sum\limits_{i=1}^{0}\dd_{x_i}:=0$ is regarded as the zero
measure. Let $\F_{\GG_0}$ be the Borel $\si$-field on $\GG_0$ induced by
the topology of weak convergence. In particular, $\GG_0$ is metrizable to be a
Polish space (cf. \cite{AKR1, AKR2, R} and references therein for geometry and analysis on $\GG_0$).

To describe  the reaction of the particle system, we first fix the death part of its transition rate.
Since the number of particles
behaves as a $Q$-process, the rate to kill $k$ particles will be $q_{|\gg|, |\gg|-k}$, where
$\gg\in\GG_0$ is the configuration of the system and $|\gg|:=\gg(E)$.
Whenever the number of
particles to be killed is fixed,  we then simply let  each particle  die at the same rate. Therefore,
the death part of the transition rate of the system reduces to

$$q_d(\gg, A):=\sum_{k=1}^{|\gg|} q_{|\gg|, |\gg|-k}
\ff{\#\{\eta\in\GG_0^{(k)}: \gg-\eta\in A\}}{\#\{\eta\in\GG_0^{(k)}:\eta\le
\gg\}},\ \ \GG_0^{(k)}:=\{\gg\in\GG_0: |\gg|=k\},\ k\ge 0.$$
where $\#$ is the cardinality of a set and $\gg-\eta\in A$ means that $\gg\ge \eta$ and $\gg-\eta\in A$.
By convention we set $q_d(0,\cdot)=0.$

Next, we go to construct the birth part of the transition rate.
Once again, since the birth rate of the particle number is determined by $Q$,
we  only need to fix the distributions of the newly produced particles.
We shall use a sequence of measures $\{\mu^{(n)}\}$ to describe
the distribution of new  particles, where $\mu^{(n)}$ is a symmetric
probability measure on $E^n, n\ge 1$.
To this end,
 we need the following assumption:

\ \newline $(H_2)\ \ \  \mu^{(n)}$ is equivalent to
$\mu^{(m)}(\cdot\times E^{m-n})$ for any $m>n\ge 1.$

\ \newline Then the birth part of the transition rate will be
determined uniquely by letting the transition rate to be symmetric
w.r.t. the probability measure (see Remark 2.1 below)

$$\pi_{\mu,\vrr}(A):= \vrr_0 1_{B}(0) +\sum_{n=1}^\infty \vrr_n \mu^{(n)}\circ
\varphi_n^{-1}(A\cap \GG_0^{(n)}),$$ where
$$\varphi_n(x):=\sum_{i=1}^n \dd_{x_i},\ \  n\ge 1, x:=(x_1,\cdots, x_n)\in E^n.$$
It is easy to see that $\varphi_n$ is continuous and hence
measurable.

We now describe the construction of  the birth part for the transition rate.
Given $m>n$, let $h_n^{(m)}$ be a fixed version of the
density of $ \mu^{(n)}$ w.r.t. $\mu^{(m)} (\cdot\times
E^{m-n})$. Let
$\mu^{(m)}((x_1,\cdots, x_n), \cdot)$ be the regular
conditional distribution of $\mu^{(m)}$ given $x_1,\cdots,
x_n$. If the system with $n$ particles $x_1, \cdots, x_n$ gives
birth to $m-n$ new particles, then we let the distribution of the new
particles  be

$$ \mu_n^{(m)}(x_1,\cdots, x_n; \cdot):= h_n^{(m)}(x_1,\cdots, x_n)
\mu^{(m)}((x_1,\cdots, x_n), \cdot).$$ It is trivial to see
that $\mu_n^{(m)}$ is symmetric in $x_1,\cdots, x_n$ and

$$ \mu^{(m)}(\d x_1,\cdots, \d x_m)= \mu^{(n)}(\d
x_1,\cdots,\d x_n) \mu_n^{(m)}(x_1,\cdots, x_n;\d
x_{n+1},\cdots, \d x_m).$$ Since $ \mu_n^{(m)}(x_1,\cdots,
x_n; A)$ is symmetric in $x_1,\cdots, x_n$, we may and will write

$$\mu_n^{(m)}(\gg, \cdot):= \mu_n^{(m)}(x_1,\cdots, x_n;
\cdot),\ \ \text{if}\ \gg=\sum_{i=1}^n\dd_{x_i}.$$ Therefore,
the birth part of the transition rate can be written as follows:

$$q_b(\gg, A):= \sum_{k=1}^\infty q_{|\gg|,|\gg|+k}
\mu_{|\gg|}^{(|\gg|+k)}(\gg, \{x\in E^k: \gg+ \varphi_k(x)\in
A\}),\ \ \ \gg\in\GG_0, A\in \F_{\GG_0}.$$

Thus,  we define the $q$-pair for the reaction of the system by letting
$q(\gg):=q(\gg,\GG_0)=q(\gg, \GG_0\setminus\{\gg\})$ and

\beq\label{rw1.2} \beg{split} q(\gg,A):= &\sum_{k=1}^\infty
q_{|\gg|, |\gg|+k}\,\mu_{|\gg|}^{(|\gg|+k)}\Big(\Big\{x\in
E^k: \gg+\varphi_k(x)\in A\Big\}\Big)\\
&+ \sum_{k=1}^{|\gg|} q_{|\gg|, |\gg|-k} \ff{\# \{\eta\in
\GG_0^{(k)}: \gg-\eta\in A\}}{\#\{\eta\in\GG_0^{(k)}: \eta\le
\gg\}},\ \ \ \gg\in\GG_0, A\in \F_{\GG_0}.\end{split}\end{equation}
This $q$-pair is regular and symmetric w.r.t. $\pi_{\mu,\vrr}$
 (see Propositions
\ref{rwP2.1} and \ref{rwP2.2} below); that is, there exists a unique $q$-process
with transition probability  kernels satisfying

\beg{equation}\label{TK}\lim_{t\to 0}\ff{P_t(\gg, A)- \dd_{\gg}(A)} t= q(\gg, A)- q(\gg) \dd_{\gg}(A)\end{equation}
for all $\gg\in\GG_0$ and $A\in \F_{\GG_0}$ such that $\lim_{t\to 0}
\sup_{\gg\in A} (1-P_t(\gg,\{\gg\}))=0,$ and the process is reversible w.r.t. $\pi_{\mu,\vrr}$.
In particular, (\ref{TK}) holds for all $A\in \F_{\GG_0}$ satisfying $\sup_{\gg\in A}|\gg|<\infty$, see e.g. \cite[Theorem 1.5(1)]{Chen}.

Since $q(\gg,\d\eta)$ is symmetric w.r.t. $\pi_{\mu,\vrr}$, the
corresponding  quadratic form is given by

\beg{equation}\label{rwD}\beg{split}\E_R^{\GG_0}(F,G) &:= \ff 1
2\int_{\GG_0\times\GG_0} ((F(\gg)-F(\eta))(G(\gg)-G(\eta))q(\gg,
\d\eta) \prr(\d\gg)\\
&= \sum_{m>n\ge 0} \int_{\GG_0^{(n)}\times \GG_0^{(m)}}
(F(\gg)-F(\eta))(G(\gg)-G(\eta)) q_b(\gg,
\d\eta)\pi_{\mu,\vrr}(\d\gg)\\
&=\sum_{n=0}^\infty \sum_{m=n+1}^\infty \vrr_nq_{n,m}
\int_{\GG_0^{(n)}}\pi_\mu(\d\gg)
\int_{E^{m-n}}(D_xF(\gg))(D_xG(\gg))\mu_n^{(m)}(\gg, \d x)
\end{split}\end{equation}
 for all $F,G$ with $\E_R^{\GG_0}(F,F)+\E_R^{\GG_0}(G,G)<\infty,$ where
$D_{x}F(\gg):=F(\gg+\varphi_{m-n}(x))-F(\gg),  \ x\in E^{m-n}.$ To
ensure that the form is well-defined in the
$L^2(\GG_0,\pi_{\mu,\vrr})$-sense, we assume that

\ \newline $(H_3)\  \mu_n^{(m)}(\gg,\cdot)$ is absolutely
continuous w.r.t. $\  \mu^{(m-n)}$ for any $m>n\ge 0$ and any
$\gg\in\GG_0^{(n)}.$

\ \newline
Under this assumption $q_b(\gg,\cdot)$ is absolutely
continuous w.r.t. $\pi_{\mu,\vrr}$, so that $\E_R^{\GG_0}$ is well-defined
on $\D(\E_R^{\GG_0}):= \{F\in L^2(\GG_0,\prr): \E_R^{\GG_0}(F,F)<\infty\}$;
that is, $\E_R^{\GG_0}(F,G)=\E_R^{\GG_0}(F',G')$ if $F,G$ represent the same classes as $F',G'$
respectively in
$L^2(\GG_0,\prr)$.  Thus, $(\E_R^{\GG_0}, \D(\E_R^{\GG_0}))$ is a
conservative symmetric Dirichlet form on $L^2(\GG_0,\prr)$
associated to the unique reversible $q$-process (see Proposition
\ref{rwP2.2} below).

If, in particular, $\mu^{(n)}:={\mu^{(1)}}^n$ for all $n\ge 1$ and
$q_{i,j}=0$ for $|i-j|>1$, the system is called a spatial birth-death
system which goes back to \cite{Preston}, see also \cite{HS} for
the study of a class of  birth-death systems on infinite
configuration spaces. In these two references  the
existence of the associated Markov processes and the description of
reversible measures were studied. Recently, there has been increasing interest in
the study of functional inequalities for spatial birth-death systems,
see e.g. \cite{Wu1} for the modified log-Sobolev inequality of
spatial birth-death systems on Poisson spaces, \cite{Kond, Wu2} for the
Poincar\'e inequality (or spectral gap) of  spatial birth-death systems
on configuration spaces. In
this paper we first study functional inequalities for the above
constructed $q$-process (i.e. the reaction process) then pass to the reaction-diffusion setting
where the particles are allowed to move dependently on $E$, i.e. undergoing interactions between them.

We prove  that if the support of
$\mu^{(1)}$ is infinite then the Dirichlet form $\E_R^{\GG_0}$ does
not satisfy the super Poincar\'e inequality (hence the associated
semigroup is not uniformly integrable, see \cite{GW, W00b}), and
it satisfies the Poincar\'e or the weak Poincar\'e inequality if
and only if so does $\E_Q$, the Dirichlet form of the $Q$-matrix
(see Theorem \ref{rwT3.1} below):

\beg{equation*}\beg{split}&\E_Q(\r, \mathbf s):=
\sum_{n=0}^\infty\sum_{m=n+1}^\infty \vrr_n q_{n,m}
(r_n-r_m)(s_n-s_m),\\
&\r=\{r_n\}, \mathbf s=\{s_n\}\in \D(\E_Q):=\{\r\in
L^2(\Z_+;\vrr): \E_Q(\r,\r)<\infty\}.\end{split}\end{equation*}
Furthermore, one has

\beq\gap(\E_Q)\label{rwgap}\ge \gap(\E_R^{\GG_0})\ge \vrr_0
\gap(\E_Q),\end{equation} where $\gap(\cdot)$ is the spectral gap
of a conservative Dirichlet form.
The first inequality in (\ref{rwgap}) follows immediately by
taking reference functions which are constant on each $\GG_0^{(n)}\
(n\ge 0)$, while to obtain the second inequality, one has to show
that, up to a multiplicative constant,  the Dirichlet form of a function $F$
dominates the square of the $L^2$-distance between  $F$ and some
function with constant value on each $\GG_0^{(n)}$, see the proof of
Theorem \ref{rwT3.1} for details. Moreover, we present an example
to show that in general one has $\gap(\E_Q)> \gap (\E_R^{\GG_0})$
(see Example 3.1 below).

Since in general $\E_R^{\GG_0}$ does not satisfy the super Poincar\'e
inequality, to derive stronger (e.g. the log-Sobolev) inequalities
one has to enlarge the Dirichlet form. To this end, we let particles
in the system move as  Markov processes. More precisely, let
$(\E_0^{(n)},\D(\E_0^{(n)}))$ be a symmetric conservative
Dirichlet form on $L^2(E^n;\mu^{(n)})$. For any function $F$ on
$\GG_0$, let $F^{(n)}:= F\circ \varphi_n$; that is,

$$F^{(n)}(x_1,\cdots, x_n):= F\Big(\sum_{i=1}^n\dd_{x_i}\Big),\ \
\ \ n\ge 1.$$  Define

$$\E_0^{\GG_0}(F,G):= \sum_{n=1}^\infty \rr_n \E_0^{(n)}(F^{(n)}, G^{(n)})$$
with $\D(\E_0^{\GG_0}):=\{F\in L^2(\GG_0,\prr): F^{(n)}\in
\D(\E_0^{(n)}), n\ge 1, \E_0^{\GG_0}(F,F)<\infty\}.$ According to
Proposition \ref{rwP2.3} below, $(\E_0^{\GG_0}, \D(\E_0^{\GG_0}))$ is a
conservative symmetric Dirichlet form on $L^2(\prr).$ Moreover,
Proposition \ref{rwP2.4} says that $\E^{\GG_0}:=\E_0^{\GG_0}+\E_R^{\GG_0}$
with domain $\D(\E^{\GG_0}):=\D(\E_0^{\GG_0})\cap \D(\E_R^{\GG_0})$ is a
symmetric Dirichlet form on $L^2(\prr).$

Now, we consider the $\phi$-variance inequality studied in \cite{W4}
(see also \cite{LO} for a special case). This inequality interpolates the
Poincar\'e and the log-Sobolev inequalities and has the additivity
property which is in particular crucial for applications in infinite dimensions.
For any probability space $(\OO, \scr B, P)$ and any
decreasing function $\phi\in C([1,2])$ with $\phi(p)>0$ for $p\in
[1,2),$ define the $\phi$-variance by

$$V_{\phi,P}(f):= \sup_{p\in [1,2)} \ff{P(f^2)- P(|f|^p)^{2/p}}{\phi(p)},
\ \ \ f\in L^2(P).$$
When $\phi\equiv 1$ and $f\ge 0$ this quantity coincides with the variance of
$f$, and when $\phi(p)= (2-p)/p$ it reduces to
$P(f^2\log f^2)$, see e.g. \cite{LO}. Thus, the following quantity
is an extension of the spectral gap and the log-Sobolev constant:

\beq\label{rw1.3}\ll_\phi(\E^{\GG_0}) := \inf\{\E^{\GG_0}(F,F): F\in
\D(\E^{\GG_0}), V_{\phi,\prr}(F)=1\}.\end{equation}
 In particular, if $\phi\equiv 1$ then $\ll_\phi=\gap$ while if $\phi(p)=(2-p)/p$
then $\ll_\phi(\E^{\GG_0})$ coincides with the log-Sobolev constant

$${\bf L}(\E^{\GG_0}):= \inf\{\E^{\GG_0}(F,F): F\in \D(\E^{\GG_0}), \text{Ent}_{\prr}(F^2)=1\}.$$
Let $\ll_\phi(\E_0^{(n)})$ and $\ll_\phi(\E_Q)$ be the
corresponding quantities of $\E_0^{(n)}$ (w.r.t. $\mu^{(n)})$ and
$\E_Q$ (w.r.t. $\vrr$). By Theorem \ref{rwT4.1} below, we have

$$\ll_\phi (\E_Q)\ge \ll_\phi(\E^{\GG_0})\ge \ll_\phi(\E_Q)\land \inf_{n\ge 1}\ll_\phi(\E_0^{(n)}).$$

Moreover, let $\E_0$ be the Dirichlet form of an independent
system on $L^2(E^\N;\mu)$; that is, $\mu={\mu^{(1)}}^\N$ and
$\E_0$ is the  sum of single Dirichlet forms
$(\E_0^{(1)}, \D(\E_0^{(1)}))$ on $L^2(E;\mu^{(1)})$.   According
to Theorem \ref{rwT4.2} below, under a mild condition, $\E^{\GG_0}$
satisfies the log-Sobolev (or the super log-Sobolev) inequality if
and only if so do $\E_0^{(1)}$ (on $L^2(E;\mu^{(1)})$) and $\E_Q$
(on $L^2(\Z_+;\vrr))$. Finally, some specific models are presented
in Section 5 to illustrate the main results.

\section{The Dirichlet Forms}

\beg{prp}\label{rwP2.1} Assume $(H_2)$. Then $(q, q(\cdot, \d
\gg))$ is a totally stable and conservative $q$-pair on $(\GG_0,\F_{\GG_0})$.\,If $(H_1)$ holds
then $q(\cdot,\d\gg)$ is symmetric w.r.t.\,$\prr$,\,i.e.\,$\prr(\d
\gg)q(\gg,\d\eta) = \prr(\d\eta) q(\eta,\d\gg).$
\end{prp}

\beg{proof} Since $Q$ is totally stable, by (\ref{rw1.2}) we have

$$q(\gg):= q(\gg,E) =\sum_{k=1}^\infty q_{|\gg|, |\gg| +k} +\sum_{k=1}^{|\gg|}
q_{|\gg|, |\gg|-k}=q_{|\gg|} <\infty,\ \ \gg\in\GG_0,$$ so that the
$q$-pair is regular too. It remains to prove the symmetry of the
measure $J(d\gg,d\eta):=\pi_{\mu,\varrho}(\d\gg)q(\gg,d\eta).$ For
any $m>n$ and  measurable sets $A_n\subset\GG_0^{(n)},
B_m\subset\GG_0^{(m)}$, let $\tt A_n:=\varphi_n^{-1}(A_n)$ and $\tt
B_m:= \varphi_m^{-1}(B_m)$. By (\ref{rw1.2}) we have

\beg{equation}\label{rw2.1}\beg{split} J(A_n\times
B_m)=&\int_{A_n} q_b(\gg, B_m)
\prr(\d\gg)\\
=&q_{n,m}\int_{A_n}\mu_n^{(m)}(y;\{x\in E^{m-n}:
\gg+\varphi(x)\in
B_m\})\prr(\d\gg)\\
=& q_{n,m}\vrr_n\int_{\tt{A}_n}\mu_n^{(m)}(y; \{x\in E^{m-n}:
(y, x)\in\tt{B}_m\})
\mu^{(n)}(\d y)\\
=&q_{n,m}\vrr_n\mu^{(m)}((\tt{A}_n\times E^{m-n})\cap\tt{B}_m).
\end{split}\end{equation}
On the other hand,

\beg{equation*} \beg{split} &J(B_m\times
A_n)=\int_{B_m}q_d(\gg,A_n)\pi_{\mu,\varrho}(\d\gg)\\
&=q_{m,n}\vrr_m\int_{\tt{B}_m} \ff{1}{\#\{\eta\in\GG_0^{(n)}:
\eta\le \varphi_m(x)\}}\sum_{\eta\in\GG_0^{(n)}: \eta \le
\varphi_m(x)} 1_{{A}_n}(\eta)\mu^{(m)}(\d x).
\end{split}\end{equation*}
Since $\tt A_n, \tt B_m$ and $\mu^{(m)}$ are symmetric in
coordinates, one has

\beg{equation*}\beg{split} J(B_m\times A_n) = &q_{m,n}\vrr_m
\int_{\tt B_m} 1_{\tt A_n}
(x_1,\cdots, x_n)\mu^{(m)}(\d x_1,\cdots,\d x_m)\\
=&q_{m,n}\vrr_m \mu^{(m)}((\tt A_n\times E^{m-n})\cap \tt
B_m).
\end{split}\end{equation*}
Combining this with (\ref{rw2.1}) and $(H_1)$, we obtain
$J(A_n\times B_m)=J(B_m\times A_n)$. Therefore, for any measurable
sets $A$ and $B$, letting $A_n:=A\cap \GG_0^{(n)}$ and $B_m:=B\cap
\GG_0^{(m)}$, we have

$$
 J(A\times B)=\sum_{n,m=0}^{\infty}J(A_n\times B_m)
 =\sum_{n,m=0}^{\infty}J(B_m\times A_n)=J(B\times A).$$\end{proof}

\paragraph{Remark 2.1.} Once the death part $q_d(\gg,\d\eta)$ is given,
$q(\gg, \d\eta)$ (hence the birth part $q_b(\gg,\d\eta)$) is
uniquely determined by its symmetry w.r.t. $\prr$. Indeed, since
the measure $J(\d\gg,\d\eta):=
\pi_{\mu,\varrho}(\d\gg)q(\gg,\d\eta)$ is symmetric, it is
uniquely determined by $J(A_n\times B_m)$ for $n>m$ and $A,B\in
\F_{\GG_0}$, which depends only on $q_d$ and $\prr$. Then $q(\gg,\d\eta)/q(\gg)$,
as the regular conditional measure of the probability measure
$J(\d\gg,\d\eta)/q(\gg)$ given $\GG_0$, is uniquely determined by
$q_d$ for $\prr$-a.e. $\GG_0$ (note that $\prr(\d\gg)$ is the first
marginal measure of $J(\d\gg,\d\eta)/q(\gg)$).

\beg{prp}\label{rwP2.2} Assume $(H_1)-(H_3)$. Then there exists a
unique $q$-process, which is reversible w.r.t. $\pi_{\mu,\vrr}$
with Dirichlet form  $(\E_R^{\GG_0}, \D(\E^{\GG_0}_R))$.
\end{prp}

\beg{proof} According to \cite[Lemma 6.52]{Chenbook}, $(\E_R^{\GG_0},
\D(\E_R^{\GG_0}))$ is a Dirichlet form on $L^2(\GG_0, \prr).$ Since
$q(\gg)= q(\gg, E)$ by \cite[Theorem 3.8]{Chenbook}, to prove the
uniqueness of the $q$-process, it suffices to verify that for any
bounded nonnegative measurable function $F$ such that

\beq\label{rwL}L_R^{\GG_0} F(\gg):= \int_{\GG_0}
(F(\eta)-F(\gg))q(\gg,\d\eta)=F,\end{equation} one has $F=0$, or
equivalently, $\text{dim}\scr U_\ll =0$ for all $\ll>0$ as
indicated by \cite[Theorem 2.37]{Chenbook}. To this end, let
$r_0:= F(0)$ and $r_n:= \ff 1 {\vrr_n} \int_{\GG_0^{(n)}} F\d \prr$
for $n\ge 1.$ We first prove that $r_n=0$ for all $n\ge 0.$ By
(\ref{rwL}) and the symmetry of the $q$-pair, we have

\beg{equation*}\beg{split} r_n&:= \int_{\GG_0^{(n)}\times \GG_0}
(F(\eta)-F(\gg))\prr(\d\gg) q(\gg, \d\eta)\\
&= \sum_{m=n+1}^\infty \vrr_n q_{n,m} (r_m-r_n)+ \sum_{m=0}^{n-1}
\vrr_m q_{m,n} (r_m-r_n)\\
&= \sum_{m=0}^\infty\vrr_n q_{n,m}
(r_m-r_n).\end{split}\end{equation*} Since the $Q$-process is
unique, by this and \cite[Theorem 6.42]{Chenbook} we have $r_n=0$
for all $n\ge 0$, that is, $F=0\ \prr$-a.e. Next, it follows from
$(H_3)$ that $q_b(\gg, \cdot)$ is absolutely continuous w.r.t.
$\prr$. Hence (\ref{rwL}) implies that

$$\sum_{n=0}^{|\gg|-1} \ff {q_{|\gg|, n}} {\#\{\eta\le\gg:
|\eta|=n\}} \sum_{\eta\le\gg: |\eta|=n} F(\eta) =
F(\gg)(q(\gg)+1).$$ Since $\vrr_0>0$ and $F=0\ \prr$-a.e., one has
$F(0)=0$ and hence by inducion in $n$ and applying the above
formula, we prove that $F\equiv 0$. Thus, the $q$-process is
unique and according to \cite[Theorems 6.7 and 6.56]{Chenbook},
the unique $q$-process is $\prr$-reversible and associated to the
Dirichlet form $(\E_R^{\GG_0},\D(\E_R^{\GG_0})).$ \end{proof}

\beg{prp}\label{rwP2.3} $(\E_0^{\GG_0},\D(\E_0^{\GG_0}))$ is a symmetric
Dirichlet form on $L^2(\prr).$
\end{prp}

\beg{proof} Obviously, $\D(\E_0^{\GG_0})$ contains the set

$$\C_1:= \{F\in \C_0: F^{(n)}\in \D(\E_0^{(n)}),\ n\ge 1\}.$$
Since $\D(\E_0^{(n)})$ is dense in $L^2(E^n;\mu^{(n)})$, $\C_1$ is
dense in $\C_0$ and hence in $L^2(\GG_0,\prr)$. Next, the
sub-Markovian property and the symmetry follow from that of $\E_0^{(n)}
(n\ge 1)$. So, it remains to verify the closedness. Let
$\{F_k\}_{k\ge 1}$ be a Cauchy sequence w.r.t. the corresponding
Sobolev norm, and let $F$ be its limit in $L^2(\GG_0,\prr)$. By the
definition of $\E_0^{\GG_0}$ one concludes that for all $m\ge 1,\
\{F_k^{(m)}\}$ is a Cauchy sequence w.r.t. the Sobolev norm
induced by $\E_0$. Since $(\E_0,\D(\E_0))$ is a Dirichlet form, it
follows that $F^{(m)}\in\D(\E_0)$ and $F_k^{(m)}\to F^{(m)}$ in
the Sobolev norm for each $m\ge 1.$ Therefore, by Fatou's lemma we
obtain

$$\sum_{n=1}^\infty \vrr_n \E_0^{(n)}(F^{(n)}, F^{(n)})
=\sum_{n=1}^\infty \liminf_{k\to\infty}\vrr_n
\E_0^{(n)}(F^{(n)}_k,F^{(n)}_k)\le \liminf_{k\to\infty}
\E_0^{\GG_0}(F_k,F_k)<\infty.$$ Hence $F\in \D(\E_0^{\GG_0})$. By using
Fatou's lemma again, we obtain $\E_0^{\GG_0}(F-F_k,F-F_k)\to 0$
as $k\to\infty.$
\end{proof}

\beg{prp}\label{rwP2.4} Let
$\D(\E^{\GG_0}):=\D(\E_0^{\GG_0})\cap\D(\E_R^{\GG_0})$ and
$\E^{\GG_0}(F,G):=\E_0^{\GG_0}(F,G)+\E_R^{\GG_0}(F,G)$. Then $(\E^{\GG_0},
\D(\E^{\GG_0}))$ is a symmetric Dirichlet form on $L^2(\prr)$.
\end{prp}

\beg{proof} Simply note that $\D(\E^{\GG_0})\supset \C_1$ and hence is
dense in $L^2(\GG_0,\prr)$. For the proof of closedness  see  \cite[Chapter I, Section 3]{MR}.
\end{proof}

\section{Poincar\'e and weak Poincar\'e inequalities for the reaction process}

We first consider the spectral gap of $(\E^{\GG_0}_R,\D(\E^{\GG_0}_R))$:
$$
\gap(\E^{\GG_0}_R):=\inf\big\{\E_R^{\GG_0}(F,F):F\in\D(\E_R^{\GG_0}),\prr(F^2)=1,\prr(F)=0\big\}.
$$
Since $\E_R^{\GG_0}$ is induced by the $Q$-matrix, it is natural for
us to relate its spectral gap  to that of $\E_Q$:

$$\gap(\E_Q):=\inf\big\{\E_Q(\r,\r):\r \in\D(\E_Q),\vrr(\r)=0,\vrr(\r^2)=1\big\}.
$$

Next, we consider the weak Poincar\'e inequality introduced in \cite{RW1}, which describes the general convergence rate
of the associated semigroup:

\beq\label{rwWPJ} \prr(F^2)\le \aa_R(r) \E_R^{\GG_0}(F,F)+
r\|F\|_\infty^2,\ \ \ r>0,\ F\in\D(\E_R^{\GG_0}), \
\prr(F)=0,\end{equation} where $\aa_R: (0,\infty)\to (0,\infty)$
is a positive function. Similarly,  this inequality is related to
the corresponding one for $\E_Q$:

\beq\label{rwWPQ} \vrr(\r^2)\le \aa_Q(r) \E_Q(\r,\r)+ r\sup_{n\ge
0} r_n^2,\ \ \ r>0,\ \r=\{r_n\}_{n\ge 0}\in\D(\E_Q), \
\vrr(\r)=0.\end{equation}

\beg{thm} \label{rwT3.1} Assume $(H_1)-(H_3)$.

$(1) \ \gap(\E_Q)\geq\gap(\E_R^{\GG_0})\geq \vrr_0 \gap(\E_Q).$
Consequently, $\gap(\E_R^{\GG_0})>0$ if and only if $\gap(\E_Q)>0.$ In
particular, for the birth-death case where $l_n:= q_{n,n+1}>0$ for
all $n\ge 0$ but $q_{n,m}=0$ for $m>n+1$, one has
$\gap(\E_R^{\GG_0})>0$ if and only if

\beq\label{rw3.0} \sup_{n\ge 0}\vrr([n+1,\infty)) \sum_{j=0}^n \ff
1 {\vrr_jl_j}<\infty.\end{equation}

$(2)\ \E_R^{\GG_0}$ satisfies the weak Poincar\'e inequality if and
only if so does $\E_Q$. More precisely, $(\ref{rwWPJ})$ implies
$(\ref{rwWPQ})$ for $\aa_Q=\aa_R$ while $(\ref{rwWPQ})$ implies
$(\ref{rwWPJ})$ for $\aa_R(r)= \ff 1{\vrr_0} \aa_Q((\vrr_0r)/4).$

$(3)$ If the support of $\mu^{(1)}$ is infinite, then $\E_R^{\GG_0}$
does not satisfy the super Poincar\'e inequality, i.e. the
following inequality does not hold for any $\bb: (0,\infty)\to
(0,\infty)$:

\beq\label{rwSPJ} \prr(F^2)\le r\E_R^{\GG_0}(F,F) +
\bb(r)\prr(|F|)^2,\ \ \ r>0, F\in \D(\E_R^{\GG_0}).\end{equation}
\end{thm}

\beg{proof} (1) For any $\r=\{r_n\}$ with $\vrr(\r)=0$ and
$\vrr(\r^2)=1,$ let $F:= \sum_{n=0}^\infty r_n1_{\GG_0^{(n)}}$. We
have $\prr(F)=0$ and $\prr(F^2)=1$, and by (\ref{rwD}),

$$\E_R^{\GG_0}(F,F) =\sum_{n=0}^\infty\sum_{m=n+1}^\infty \vrr_nq_{n,m} (r_n-r_m)^2
=\E_Q(\r,\r).$$
Then

$$\gap(\E_R^{\GG_0})\le \inf\{\E_Q(\r,\r): \vrr(\r)=0, \vrr(\r^2)=1\}=\gap(\E_Q).$$
Next, by the triangle inequality of the $L^2$-norm,

\beq\label{rw3.1} \beg{split}\E_R^{\GG_0}(F,F)
 &=\sum_{n=0}^\infty\sum_{m=n+1}^\infty \vrr_nq_{n,m}
\int_{E^m} (F^{(m)}(x,y) - F^{(n)}(y))^2 \mu_n^{(m)}(y;\d x)\mu^{(n)}(\d y)\\
&\ge \sum_{n=0}^\infty\sum_{m=n+1}^\infty \vrr_nq_{n,m}
\Big(\ss{\mu^{(m)}({F^{(m)}}^2)}
-\ss{\mu^{(n)}({F^{(n)}}^2)}\Big)^2.\end{split}\end{equation}
Let

$$\ll_0(\E_Q):= \inf\{\E_Q(\r,\r): \vrr(\r^2)=1, r_0=0\}.$$
It follows from (\ref{rw3.1}) that if $F(0)=0$ then

$$\E_R^{\GG_0}(F,F)\ge \ll_0(\E_Q) \sum_{n=0}^\infty \vrr_n \mu^{(n)}({F^{(n)}}^2)=
\ll_0(\E_Q) \prr(F^2).$$ Thus, for any $F\in \D(\E_R^{\GG_0})$ with
$\prr(F)=0$,

$$\E_R^{\GG_0}(F,F)= \E_R^{\GG_0}(F-F(0), F-F(0)) \ge \ll_0(\E_Q)\prr((F-F(0))^2)
\ge \ll_0(\E_Q)\prr(F^2).$$ This implies that $\gap(\E_R^{\GG_0})\ge
\ll_0(\E_Q).$ Since for any $\r$ with $r_0=0$ and $m_\vrr(\r^2)=1$
one has

$$\vrr(\r^2)-\vrr(\r)^2 \ge \vrr(\r^2) - \vrr(\r^2) (1-\vrr_0)= \vrr_0,$$
$\ll_0(\E_Q)\ge \vrr_0\gap(\E_Q)$ and hence the desired lower
bound of $\gap(\E_R^{\GG_0})$ follows. Therefore, the proof of (1) is
finished by noting that  for the birth-death case one has
$\gap(\E_Q)>0$ if and only if (\ref{rw3.0}) holds, see
\cite{Miclo} or \cite{Chen}.

(2) By taking $F:=\sum_{n=0}^\infty r_n1_{\GG_0^{(n)}}$ one
concludes that (\ref{rwWPJ}) implies (\ref{rwWPQ}) for
$\aa_Q=\aa_R$. On the other hand,  for any $F\in \D(\E_R^{\GG_0})$
with $F(0)=0$, it follows from (\ref{rw3.1}) and (\ref{rwWPQ})
that

\beg{equation*}\beg{split}\aa_Q(r) \E_R^{\GG_0}(F,F) &\ge \prr(F^2)
-\Big(\sum_{n=1}^\infty \vrr_n
\ss{\mu^{(n)}({F^{(n)}}^2)}\Big)^2
-r \|F\|_\infty^2\\
&\ge \prr(F^2) \vrr_0 -r\|F\|_\infty^2.\end{split}\end{equation*}
Therefore, for any $F\in \D(\E_R^{\GG_0})$ with $\prr(F)=0$,

\beg{equation*}\beg{split} \prr(F^2)&\le \prr((F-F(0))^2) \le \ff
{1}{\vrr_0}
\aa_Q(r) \E_R^{\GG_0}(F,F) +\ff {r}{\vrr_0}  \|F-F(0)\|_\infty^2\\
&\le \ff{1}{\vrr_0}\aa_Q(r) \E_R^{\GG_0}(F,F) +\ff{4 r}{\vrr_0}
\|F\|_\infty^2.\end{split}
\end{equation*}
This implies (\ref{rwWPJ}) for $\aa_R=
\ff{1}{\vrr_0}\aa_Q(\vrr_0r/4).$

(3) For any nonnegative $f\in L^2(\mu)$, let
$F(\gg):=\gg(f)1_{\GG_0^{(1)}}$. Then $$\prr (F)=\vrr_1
\mu^{(1)}(f), \prr(F^2)=\vrr_1 \mu^{(1)}(f^2)$$ and

$$ \E_R^{\GG_0}(F,F) =\vrr_0 q_{0,1} \mu^{(1)}(f^2) + \vrr_1 \sum_{m=2}^\infty q_{1,m}
\mu^{(1)}(f^2) \le (\vrr_0\lor\vrr_1)q_1 \mu^{(1)}(f^2).$$ Thus,
if the super Poincar\'e inequality holds then there exists a
constant $c>0$ such that $\mu^{(1)}(f^2)\le c\mu^{(1)}(f)^2$ for
all nonnegative $f$, which is impossible if the support of
$\mu^{(1)}$ is infinite.\end{proof}

\paragraph{Remark 3.1.} Let $\pi_\si$ be the Poisson measure with ($\si$-finite)
intensity $\si$. It is well-known that the following Poincar\'e
inequality holds (see  \cite[Remark 1.4]{Wu1}):

\beq\label{rwwu}\pi_\si (F^2) \le \int_{\GG_0} \d\pi_\si \int_E(D_x
F)^2 \si(\d x)+ \pi_\si (F)^2,\ \ \ F\in
L^2(\pi_\si).\end{equation} See  \cite{BCC, Kond, Wu2} for
extensions to a class of Gibbs measures with $E=\R^d$. Thus, in
our present setting one has $\gap(\E_R^{\GG_0})\ge 1$ provided
$\vrr\equiv 1$ and $l_n:= q_{n,n+1}\ge 1, q_{n,m}=0$ for $m>n+1$.
But it is easy to see that in this case (\ref{rw3.0}) holds if and
only if  $\inf_{n\ge 0}nl_n>0$. Therefore, Theorem \ref{rwT3.1}
(1) provides a much weaker and sharp condition for
$\gap(\E_R^{\GG_0})>0$.






To conclude this section, let us present an example to show that in general $\gap(\E_R^{\GG_0})$ is
strictly less than $\gap(\E_Q)$.

\ \newline {\bf Example 3.1.} Let $q_{1,k}= \bb_k>0$ and
$q_{k,1}=\ff 1 2$ for $k\ne 1,$ and $q_{i,j}=0$ for $i,j\ne 1.$ By
$(H_1)$ one has $\vrr_1= (1+2q_1)^{-1}$ and $\vrr_k= 2\vrr_1 \bb_k
(k\ne 1)$, where $q_1:= \sum_{k\ne 1} q_{1,k}<\infty.$ Then
$\gap(\E_Q)=\ff 1 2$ (see \cite[Example 4.7]{CW}). On the other
hand, if $\mu^{(1)}$ is non-trivial, then there exists  $f\in
L^2(E;\mu^{(1)})$ with $\mu^{(1)}(f)=0$ and $\mu^{(1)}(f^2)=1$.
Let $F(\gg):= \gg(f)1_{\GG_0^{(1)}}(\gg).$ We have $\prr(F)=0$ and
$\prr(F^2)= \vrr_1$. Moreover, by the symmetry of the $q$-pair,

$$\E_R^{\GG_0}(F,F)= \sum_{k\ne 1} q_{k,1}\vrr_k = \ff 1 2(1-\vrr_1).$$
Therefore,

$$\gap(\E_R^{\GG_0})\le \ff{1-\vrr_1}{2\vrr_1}<\ff 1 2=\gap(\E_Q),\ \
\text{if}\ \vrr_1>\ff 1 2.$$

\section{Functional inequalities for the reaction-diffusion\\ process}

We first consider $\ll_\phi(\E^{\GG_0})$ which provides a certain
exponential convergence rate of the corresponding Markov
semigroup, see \cite{W4}.

\beg{thm}\label{rwT4.1} Assume $(H_1)-(H_3)$. Let
$\ll_\phi(\cdot)$ be the quantity  defined as in $(\ref{rw1.3})$
for a Dirichlet form. We have
$$\ll_\phi(\E_Q)\ge \ll_\phi(\E^{\GG_0})\ge \ll_\phi(\E_Q)\land
 \inf_{n\ge 1}\ll_\phi(\E_0^{(n)}).$$ Consequently,
if $ \inf_{n\ge 1}\ll_\phi(\E_0^{(n)}) \ge \ll_\phi(\E_Q)$, then
$\ll_\phi(\E^{\GG_0})=\ll_\phi(\E_Q).$\end{thm}

\beg{proof} Let $F\in \D(\E^{\GG_0})$ with $V_{\phi,\prr}(F)=1.$ We
have

\beq\label{rw4.1}\beg{split} \E_0^{\GG_0}(F,F)
&:=\sum_{n=1}^\infty \vrr_n \E_0^{(n)}(F^{(n)}, F^{(n)}) \ge
 \inf_{n\ge 1}\ll_\phi(\E_0^{(n)}) \sum_{n=1}^\infty
\vrr_n V_{\phi,\prr} (F^{(n)})\\
&\ge \sup_{p\in [1,2)} \ff{ \inf_{n\ge
1}\ll_\phi(\E_0^{(n)})}{\phi(p)} \Big(\prr(F^2) -\vrr_0F(0)^2
-\sum_{n=1}^\infty \vrr_n
\mu^{(n)}(|F^{(n)}|^p)^{2/p}\Big).\end{split}\end{equation} Next,
letting $r_0:= F(0)$ and $r_n:= \mu^{(n)}(|F^{(n)}|^p)^{1/p}$ for
$n\ge 1$, we have

\beq\label{rw4.2}\beg{split} &\vrr_0 F(0)^2 + \sum_{n=1}^\infty
\vrr_n
\mu^{(n)}(|F^{(n)}|^p)^{2/p} = \sum_{n=0}^\infty \vrr_n r_n^2\\
&\le \Big(\sum_{n=0}^\infty \vrr_nr_n^p\Big)^{2/p}
+\ff{\phi(p)\E_Q(\r,
\r)}{\ll_\phi(\E_Q)}\\
&= \prr(|F|^p)^{2/p} +\ff {\phi(p)}{\ll_\phi(\E_Q)}
\sum_{n=0}^\infty\sum_{m=n+1}^\infty \vrr_n q_{n,m} (r_m-r_n)^2.
\end{split}\end{equation}
Since $F^{(k)}$ is symmetric, we have $\mu^{(k)}(|F^{(k)}|^p) =
\  \mu^{(k)}(|F^{(k)}|^p)$ for any $k\ge 1.$ Then, by the
triangle inequality and  Jensen's inequality,

\beg{equation*}\beg{split}(r_n-r_m)^2
&=\big(\mu^{(n)}(|F^{(n)}|^p)^{1/p}
-\mu^{(m)}(|F^{(m)}|^p)^{1/p}\big)^2\\
& \le \bigg(\int_{E^m}
|F^{(m)}(x,y)- F^{(n)}(y)|^p \mu_n^{(m)}(y;\d x)\mu^{(n)}(\d y)\bigg)^{2/p}\\
&\le \int_{E^m} (F^{(m)}(x,y)-F^{(n)}(y))^2 \mu_n^{(m)}(y;\d
x)\mu^{(n)}(\d y).\end{split}\end{equation*} Thus,

$$\sum_{n=0}^\infty\sum_{m=n+1}^\infty \vrr_n q_{n,m} (r_m-r_n)^2
\le \E_R^{\GG_0}(F,F).$$ Combining this with (\ref{rw4.1}) and
(\ref{rw4.2}), we obtain

$$\E_0^{\GG_0}(F,F) \ge  \inf_{n\ge 1}\ll_\phi(\E_0^{(n)}) \bigg\{ V_{\phi,\prr}(F) -\ff
{\E_R^{\GG_0}(F,F)}{\ll_\phi(\E_Q)}\bigg\}.$$ Equivalently,

$$V_{\phi,\prr}(F) \le \ff{\E_0^{\GG_0}(F,F)}{ \inf_{n\ge 1}\ll_\phi(\E_0^{(n)})}+\ff{\E_R^{\GG_0}(F,F)}
{\ll_\phi(\E_Q)} \le \ff{\E^{\GG_0}(F,F)}{ \inf_{n\ge
1}\ll_\phi(\E_0^{(n)})\land \ll_\phi(\E_Q)}.$$ This implies that
$\ll_\phi(\E^{\GG_0})\ge  \inf_{n\ge
1}\ll_\phi(\E_0^{(n)})\land\ll_\phi(\E_Q).$ Finally, for any
$\r\in \D(\E_Q)$, let $F:= \sum_{n=0}^\infty r_n 1_{\GG_0^{(n)}}.$
We have $\E_0^{\GG_0}(F,F)=0$ since $F^{(n)}:= r_n, n\ge 1.$
Moreover, $V_{\phi,\prr}(F)=V_{\phi,\vrr}(\r).$ Hence
$\ll_\phi(\E_Q)\ge \ll_\phi (\E^{\GG_0}).$ \end{proof}

Obviously, if  $\phi\in C[1,2]$ is strictly decreasing with
$\phi(2)=0$ then $\ll_\phi(\E^{\GG_0})>0$ implies the following super
Poincar\'e inequality for some positive function $\bb:$

\beq\label{rwSP} \prr(F^2)\le r \E^{\GG_0}(F,F)+ \bb(r)
\prr(|f|)^2,\ \ \ \ r>0, F\in \D(\E^{\GG_0}).\end{equation} Thus,
according to Theorem \ref{rwT3.1}, if $\mu^{(1)}$ is not finitely
supported, then the non-triviality of $\E_0$ is  necessary for
$\ll_\phi(\E^{\GG_0})>0.$ But in general, $\ll_\phi(\E^{\GG_0})>0$
only implies a certain functional inequality of $\E_0$ rather than
$ \inf_{n\ge 1}\ll_\phi(\E_0^{(n)})>0$. To see this, let us
consider a simple situation where $\E_0$ is the Dirchlet form of
an independent particle system. More precisely, let
$\mu={\mu^{(1)}}^\N$ be the product measure and let $(\E_0^{(1)},
\D(\E_0^{(1)}))$ be a conservative symmetric Dirichlet form on
$L^2(E;\mu^{(1)})$. For any $f\in L^2(E^\N;\mu)$ and $x=
(x_1,x_2,\cdots)\in E^\N$, let

$$f_{x|i}(y):= f(x_1,\cdots, x_{i-1}, y, x_{i+1},\cdots),\ \ \
i\ge 1, y\in E.$$ Define

\beg{equation}\label{rw4.3}\beg{split}&\E_0(f,g) :=
\sum_{i=1}^\infty
 \int_{E^\N} \E_0^{(1)}(f_{x|i},g_{x|i})\mu(\d x),\\
&\D(\E_0):= \{f\in L^2(\mu): f_{x|i}\in \D(\E_0^{(1)}),\
\mu\text{-a.e.}\ x, i\ge 1, \E_0^{(n)}(f,f)<\infty\}.
\end{split}\end{equation}
Then it is easy to see that $(\E_0,\D(\E_0))$ is a symmetric
conservative Dirichlet form on $L^2(E^\N;\mu).$ Moreover, since a
function on $E^n$ can be regarded as a cylindrical function on
$E^\N$, we have the following Dirichlet forms:

\beq\label{rwN} \E_0^{(n)}(f,g):= \E_0(f,g),\  \D(\E_0^{(n)}):=
\{f\in L^2(E^n;(\mu^{(1)})^n): \E_0(f,f)<\infty\}, \ \ \ \ n\ge
1.\end{equation}

We study the log-Sobolev inequality

\beq\label{rwLS} \prr(F^2\log F^2)\le C_1(\E^{\GG_0})\E^{\GG_0}(F,F) +
C_2(\E^{\GG_0}),\ \ \ F\in \D(\E^{\GG_0}), \prr(F^2)=1\end{equation} by
using the following corresponding ones:

\beq \label{rwLS0} \mu^{(1)}(f^2\log f^2)\le
C_1(\E_0^{(1)})\E_0^{(1)}(f,f) +C_2(\E_0^{(1)}),\ \ \ f\in
\D(\E_0^{(1)}), \mu^{(1)}(f^2)=1,\end{equation} \beq\label{rwLSQ}
\vrr(\r^2\log\r^2) \le C_1(\E_Q) \E_Q(\r,\r) +C_2(\E_Q),\ \ \
\r\in \D(\E_Q),\ \vrr(\r^2)=1.\end{equation}

\beg{thm}\label{rwT4.2} Assume $(H_1)-(H_3)$ and let $\E_0^{(n)}$
be given by $(\ref{rw4.3})$ and $(\ref{rwN})$.

$(1)$ Assume that $(\ref{rwLS0})$ and $(\ref{rwLSQ})$  hold. If
there exists $\dd>0$ such that $\sum_{n=0}^\infty \vrr_n \e^{\dd
n}<\infty$ then  $(\ref{rwLS})$ holds for

\beg{equation*}\beg{split} &C_1(\E^{\GG_0}) = C_1(\E_0)
\lor \big\{\big(1+ \dd^{-1}C_2(\E_0)\big)C_1(\E_Q)\big\},\\
& C_2(\E^{\GG_0})= C_2(\E_Q) (1+\dd^{-1} C_2(\E_0))+
\dd^{-1}C_2(\E_0)\sum_{n=0}^\infty \vrr_n \e^{\dd n-1}.
\end{split}\end{equation*}

$(2)$ If  $(\ref{rwLS})$ holds  then $(\ref{rwLSQ})$ holds for
$C_1(\E_Q)=C_1(\E^{\GG_0})$ and $C_2(\E_Q)= C_2(\E^{\GG_0})$, and
$(\ref{rwLS0})$ holds for $C_1(\E_0^{(1)})=C_1(\E^{\GG_0})$ and

$$C_2(\E_0^{(1)})= \inf_{n\ge 1}\ff 1 n \big\{\log \vrr_n +
C_1(\E^{\GG_0})q_n + C_2(\E^{\GG_0})\big\}.$$
\end{thm}

\beg{proof} (1)  Let $F\in \D(\E^{\GG_0})$ with $\prr(F^2)=1$. By
the sub-additivity property of the entropy (see e.g.
\cite[(4.2)]{Ledoux}), for any $F\in \D(\E_{\GG_0})$ we have

$$\text{Ent}_{\mu^{(n)}}({F^{(n)}}^2)
\leq
\sum_{i=1}^{n}\int_{E^{n}}\text{Ent}_{\mu^{(1)}}({F^{(n)}_{x|i}}^2)\mu^{(n)}(\d
x).$$ Then by (\ref{rwLS0}) we have

$$\text{Ent}_{\mu^{(n)}}({F^{(n)}}^2)
\le C_1(\E_0)\E_0 (F^{(n)}, F^{(n)}) + n C_2(\E_0)
\mu^{(n)}({F^{(n)}}^2).$$ Thus,

\beq\label{rw4.6} \beg{split}\prr(F^2\log F^2) &\le
\sum_{n=0}^\infty \vrr_n \mu^{(n)}({F^{(n)}}^2)
\log\mu^{(n)}({F^{(n)}}^2) \\
&\quad + C_1(\E_0) \E_0^{\GG_0}(F,F)+C_2(\E_0)\sum_{n=0}^\infty
n\vrr_n \mu^{(n)}(
{F^{(n)}}^2)\\
&\le (1+\dd^{-1} C_2(\E_0))\sum_{n=0}^\infty \vrr_n
\mu^{(n)}({F^{(n)}}^2)
\log \mu^{(n)}({F^{(n)}}^2)\\
&\quad +C_1(\E_0)\E_0^{\GG_0}(F,F) + \dd^{-1}C_2(\E_0)
\sum_{n=0}^\infty \vrr_n \e^{\dd n-1},\end{split}\end{equation}
where the last step is due to Young's inequality. Next, by
(\ref{rwLSQ}) we have

\beg{equation}\label{rw4.7}\beg{split}
&\sum_{n=0}^{\infty}\vrr_n\mu^{(n)}({F^{(n)}}^2)\log
\mu^{(n)}({F^{(n)}}^2)
- \prr(F^2)\log\prr(F^2)\\
&\le C_1(\E_Q)\sum_{n=0}^{\infty}\sum_{m=n+1}^\infty
\vrr_nq_{n,m} \Big(\ss{\mu^{(n)}({F^{(n)}}^2)}-\ss{\mu^{(m)}({F^{(m)}}^2)}\Big)^2\\
&\qquad+ C_2(\E_Q) \prr(F^2).
\end{split}\end{equation} We may regard $F^{(n)}$ as a function in
$L^2(E^m;\mu^{(m)})$ so that the triangle inequality and the
symmetry of $F^{(n)}$ and $F^{(m)}$ imply

\beg{equation*}\beg{split}
\Big(\ss{\mu^{(n)}({F^{(n)}}^2)}-\ss{\mu^{(m)}({F^{(m)}}^2)}\Big)^2
&=\Big(\ss{\mu^{(n)}({F^{(n)}}^2)}-\ss{\mu^{(m)}({F^{(m)}}^2)}\Big)^2\\
&\le \mu^{(m)}((F^{(n)}-F^{(m)})^2).\end{split}\end{equation*}
Therefore, combining (\ref{rw4.6}) and (\ref{rw4.7}) and noting
that $\prr(F^2)=1$, we arrive at

\beg{equation*}\beg{split} \prr(F^2\log F^2)\le & C_1(\E_0) \E_0^{\GG_0} (F,F) + (1+ \dd^{-1}C_2(\E_0))
C_1(\E_Q)\E_R^{\GG_0}(F,F)\\
&+(1+\dd^{-1}C_2(\E_0))C_2(\E_Q) +\dd^{-1}
C_2(\E_0)\sum_{n=0}^\infty \vrr_n\e^{\dd
n-1}.\end{split}\end{equation*}

(2) Assume that (\ref{rwLS}) holds. For any $\r\in \D(\E_Q)$,
letting $F:=\sum_{n=0}^\infty r_n1_{\GG_0^{(n)}}$ we have
$\E_0^{\GG_0}(F,F)=0$ and $\E_R^{\GG_0}(F,F) = \E_Q(\r,\r).$
Moreover, Ent$_{\prr} (F^2) =\text{Ent}_{\vrr}(\r^2)$ and
$\prr(F^2)=\vrr(\r^2)$. Then we obtain (\ref{rwLSQ}) for
$C_i(\E_Q)=C_i(\E_R^{\GG_0}), i=1,2.$ Next, for any $f\in
\D(\E_0^{(1)})$ with $\mu^{(1)}(f^2)=1$ and any $n\ge 1,$ let

\beq\label{rw4.8} F(\gg):= \beg{cases} f(x_1)\cdots f(x_n), &\text{if}\ \gg=\sum_{i=1}^n \dd_{x_i}\in\GG_0^{(n)},\\
0, &\text{otherwise.}\end{cases}\end{equation}
Then it is easy to see that

\beg{equation*}\beg{split}&\prr(F^2)= \vrr_n,
\text{Ent}_{\prr}(F^2) =n\vrr_n \mu^{(1)}(f^2\log f^2) -\vrr_n\log \vrr_n,\\
& \E_0^{\GG_0}(F,F) =n \vrr_n \E_0^{(1)}(f,f),\ \ \E_R^{\GG_0}(F,F)=
\sum_{l=0}^{n-1} \vrr_l q_{l,n} + \sum_{l=n+1}^\infty
\vrr_nq_{n,l}=\vrr_nq_n,\end{split}\end{equation*} where the last
equality is due to $(H_1)$. Thus, (\ref{rwLS}) implies
(\ref{rwLS0}) for the desired constants.
\end{proof}

Now, let $\phi(p):= (2-p)/p$ so that $\ll_\phi$ coincides with the
log-Sobolev constant and let $\vrr$ satisfy $\sum_{n=1}^\infty
\e^{\dd n}\vrr_n<\infty$ for some $\dd>0.$ According to Theorem
\ref{rwT4.2}, if $\ll_\phi(\E_Q)>0$ and (\ref{rwLS0}) holds, then
(\ref{rwLS}) holds. Since $\ll_\phi(\E_Q)>0$ implies
$\gap(\E_Q)>0$, we have $\gap(\E^{\GG_0})>0$ according to Theorem
\ref{rwT3.1}. Thus, $\ll_\phi(\E^{\GG_0})>0$. On the other hand,
however,  there are a lot of examples where (\ref{rwLS0}) holds
but $\ll_\phi(\E_0^{(1)})=0$ (hence, $ \inf_{n\ge
1}\ll_\phi(\E_0^{(n)})=0$). Therefore, as claimed before,
$\ll_\phi(\E^{\GG_0})>0$ does not imply $ \inf_{n\ge
1}\ll_\phi(\E_0^{(n)})>0.$

Theorem \ref{rwT4.2} enables us to study  the super log-Sobolev
inequality

\beq\label{rwSLS}\prr(F^2\log F^2)\le r \E^{\GG_0}(F,F) +\bb(r),\ \ \
F\in \D(\E^{\GG_0}), \prr(F^2)=1,\end{equation} where $\bb:
(0,\infty)\to (0,\infty)$ is a positive function. According to
\cite{DS}, this inequality is equivalent to the supercontractivity
of $P_t^{\GG_0}: \ \|P_t^{\GG_0}\|_{L^2(\GG_0,\prr)\to
L^4(\gg,\prr)}<\infty$ for all $t>0.$ We shall study this
inequality by using the corresponding ones for $\E_Q$ and
$\E_0^{(1)}$:

\beq \label{rwSLS0} \mu^{(1)}(f^2\log f^2)\le r\E_0^{(1)}(f,f)
+\bb_0(r),\ \ \ f\in \D(\E_0^{(1)}),
\mu^{(1)}(f^2)=1,\end{equation}

\beq\label{rwSLSQ} \vrr(\r^2\log\r^2) \le r \E_Q(\r,\r)
+\bb_Q(r),\ \ \ \r\in \D(\E_Q),\ \vrr(\r^2)=1.\end{equation} The
following result is a direct consequence of Theorem \ref{rwT4.2}.

\beg{cor}\label{rwC4.3} Consider the situation of Theorem $\ref{rwT4.2}$
and assume that $\sum_{n=0}^\infty \vrr_n \e^{\dd n}<\infty$ for
some $\dd>0.$ Then $(\ref{rwSLS0})$ and $(\ref{rwSLSQ})$ imply
$(\ref{rwSLS})$ for

$$\bb(r)= \bb_Q\big(\dd/(\dd+\bb_0(r))\big) (1+\dd^{-1} \bb_0(r)) +\dd^{-1}
\bb_0(r) \sum_{n=0}^\infty \vrr_n \e^{\dd n-1}.$$ On the other
hand, $(\ref{rwSLS})$ implies $(\ref{rwSLSQ})$ for $\bb_Q=\bb$ and
$(\ref{rwSLS0})$ for

$$\bb_0(r)= \inf_{n\ge 1} \ff 1n \{\log \vrr_n
+\bb(r) +rq_n\}.$$\end{cor}

Finally, the above arguments can be also applied to the super
Poincar\'e inequality.

\beg{cor} \label{rwC4.5} Consider the situation of Theorem
$\ref{rwT4.2}$ and assume that  $\gap(\E_0)>0$. Then $\E^{\GG_0}$
satisfies $(\ref{rwSP})$ for some $\bb$ if and only if there exist
$\bb_0, \bb_Q: (0,\infty)\to (0,\infty)$ such that

\beq \label{rwSP0} \mu^{(1)}(f^2)\le r\E_0^{(1)}(f,f)
+\bb_0(r)\mu^{(1)}(|f|)^2,\ \ \ f\in
\D(\E_0^{(1)}),r>0,\end{equation}

\beq\label{rwSPQ} \vrr(\r^2) \le r \E_Q(\r,\r)
+\bb_Q(r)\vrr(|\r|)^2,\ \ \ \r\in \D(\E_Q), r>0.\end{equation}
\end{cor}

\beg{proof} The proof that (\ref{rwSP}) implies (\ref{rwSP0}) and
(\ref{rwSPQ}) is similar to the proof that (\ref{rwLS}) implies
(\ref{rwLS0}) and (\ref{rwLSQ}), so we only prove the converse.
Since the super Poincar\'e inequality is equivalent to a Sobolev
type inequality, that is, replacing the function $\log$ in the
log-Sobolev inequality by some function increasing to infinity as
the variable goes to infinity (see \cite{GW} or \cite{W00b}), and
since $\gap(\E_0)>0$, by \cite[Theorem 1.1]{W4} and (\ref{rwSP0})
we have $ \inf_{n\ge 1}\ll_\phi(\E_0^{(n)})>0$ for some strictly
decreasing $\phi\in C([1,2])$ with $\phi(2)=0$. For any $F\in
\D(\E^{\GG_0})$, by the sub-additivity of $V_{\phi,\mu}$ we have

$$V_{\phi,\mu^{(n)}}(F^{(n)})\le \ff 1 {\ll_\phi (\E_0^{(1)})} \E_0^{(n)}
(F^{(n)}, F^{(n)}),\ \ \ n\ge 1.$$
This implies

\beq\label{rw4.13} \prr(F^2)\le \sum_{n=0}^\infty \vrr_n
\mu^{(n)}(|F^{(n)}|^p)^{2/p}
+\ff{\phi(p)}{\ll_\phi(\E_0^{(1)})}\E_0^{\GG_0}(F,F), \ \ \ p\in
[1,2).\end{equation} Next, we claim that for any probability space
$(\OO,\scr B,P)$ and any function $h\in L^2(P)$ one has

\beq\label{rw4.14} P(|h|^p)^{2/p}\le \ff 1 2 P(h^2) + \ff 2 p
(\ff4 p)^{2p/(2-p)}P(|h|)^2,\ \ \ \ \ p\in [1,2).\end{equation}
Indeed, letting $P(|h|)=1$ we have

$$P(|h|^p)^{2/p}\le \ff 2 p\big(P(|h|^p1_{\{|h|>R\}})^{2/p} + R^{2}\big)
\le \ff 2 p\big(R^{-(2-p)/p}P(h^2) + R^{2}\big),\ \ R>0.$$ Taking
$R= (\ff 4 p)^{p/(2-p)}$ we prove (\ref{rw4.14}). Letting
$c(p):=\ff 2 p (\ff4 p)^{2p/(2-p)}$, by (\ref{rw4.13}),
(\ref{rw4.14}) and  (\ref{rwSPQ}) we obtain

\beg{equation*}\beg{split} \prr(F^2)&\le
\ff{2\phi(p)}{\ll_{\phi}(\E_0^{(1)})} \E_0^{\GG_0}(F,F) +
2c(p)\sum_{n=0}^\infty
\vrr_n\mu^{(n)}(|F^{(n)}|)^2\\
&\le \ff{2\phi(p)}{\ll_{\phi}(\E_0^{(1)})} \E_0^{\GG_0}(F,F)+
2c(p)r_1\E_R^{\GG_0}(F,F) +
2c(p)\bb_Q(r_1)\prr(|F|)^2\end{split}\end{equation*} for any $
p\in [1,2)$ and any $r_1>0.$ Since $\phi(2)=0$, (\ref{rwSP}) holds
with

$$\bb(r):=\inf\Big\{2c(p)\bb_Q(r_1): p\in [1,2), r_1>0 \text{\ such\ that\ } \ff{2\phi(p)}
{ \inf_{n\ge 1}\ll_\phi(\E_0^{(n)})}\lor(2c(p)r_1)\le r\Big\}$$
which is finite for any $r>0.$\end{proof}

\section{Examples}

 In this section we present three
specific models where the underlying Markov chain is the birth-death
process; that is, $Q$ and $\vrr$ satisfy $(H_1)$ with $l_n:=
q_{n,n+1}>0$ for $n\ge 0$ and $q_{n,m}=0$ for $m>n+1$. In the
first two examples $\E_0$ refers to some infinite-dimensional
diffusion on a manifold, where in the first example the diffusion
process is without interaction but the manifold is non-compact,
and in the second example the diffusion process is given by the
one-dimensional stochastic Ising model over a compact manifold.
Finally, as a supplement to Theorem \ref{rwT3.1}(3),  we show in
the last example that the pure reaction Dirichlet form $\E_R^{\GG_0}$ may
satisfy the log-Sobolev inequality if $E$ is finite.

\paragraph{Example 5.1. (with independent diffusions)}
Let $M$ be connected and noncompact with Ricci curvature bounded from below. Let
$V\in C(M)$ such that $V + c \rr^\theta$ is bounded, where $c>0,
\theta >1$ are constants, and $\rr$ is the Riemannian distance
function to a fixed point. By the volume comparison theorem (see
\cite{CGT}) one has $Z:=\int_E \e^{V(x)}\d x<\infty$, where $\d x$
is the volume measure. Let $\mu^{(1)}(\d x):= Z^{-1} \e^{V(x)}\d
x$ and $\E_0^{(1)}(f,g):= \mu^{(1)}(\<\nn f,\nn g\>)$ with
$\D(\E_0^{(1)}):= H^{2,1}(\mu^{(1)})$, the completion of
$C_0^\infty(M)$ under the Sobolev norm
$\|\cdot\|_{L^2(\mu^{(1)})}+ \E_0^{(1)}(\cdot,\cdot)^{1/2}$. Let
$\E_0$ and $\E_0^{(n)}$ be given by (\ref{rw4.3}) and (\ref{rwN}).
We have:

\ \newline {\bf (i)} $\E_R^{\GG_0}$ (and hence $\E^{\GG_0}$) always
satisfies the weak Poincar\'e inequality, and it (equivalently,
$\E^{\GG_0}$)  satisfies the Poincar\'e inequality  if and only if
(\ref{rw3.0}) holds. \newline {\bf (ii)}  (\ref{rwSP}) holds for
some $\bb$ if and only if

\beq\label{rw5.3} \lim_{n\to\infty}\vrr([n+1,\infty)) \sum_{j=0}^n
\ff 1 {\vrr_jl_j}=0.\end{equation} {\bf (iii)} Let $\phi(p):=
(2-p)^\aa$ for $\aa\in (0,1].$ Then $\ll_\phi(\E^{\GG_0})>0$  if and
only if $\theta\ge 2/(2-\aa)$ and

\beq\label{rw1.9} \sup_{n\ge 0}\vrr([n+1,\infty))\Big(\log
{\vrr([n+1,\infty))}^{-1}\Big)^\aa \sum_{j=0}^n \ff 1
{\vrr_jl_j}<\infty.\end{equation} {\bf (iv)} Assume that
$\sum_{n=0}^\infty \vrr_n \e^{\dd n}<\infty$ for some $\dd>0$.
Then $\E^{\GG_0}$ satisfies (\ref{rwSLS}) for some $\bb$ if and only
if $\theta>2$ and

\beq\label{rw1.7} \lim_{n\to\infty}\big[\vrr([n+1,\infty))\log
\vrr([n+1,\infty))^{-1}\big] \sum_{j=0}^n \ff 1
{\vrr_jl_j}=0.\end{equation}

\beg{proof} {\bf (i)} follows from Theorem \ref{rwT3.1} (1) and
(2)
 and the following facts: any reversible irreducible
countable Markov chain satisfies the weak Poincar\'e inequality
(see \cite[Theorem 3.1]{RW1} or \cite[Corollary 1.3]{W10});
$\gap(\E_0)=\gap(\E_0^{(1)})>0$ according to \cite{W11} or
\cite[Corollary 1.3]{RW2}; $\gap(\E_Q)>0$ if and only if
(\ref{rw3.0}) holds (see \cite{Miclo} or \cite{Chen}).

{\bf (ii)} follows from Corollary \ref{rwC4.5}  and the facts that
$\gap(\E_0^{(1)})>0$ and $\E_0^{(1)}$ satisfies the super
Poincar\'e inequality since $\theta>1$ (see \cite[Corollary
2.5]{W3} or \cite[Corollary 1.3]{RW2}), while by the discrete
Hardy inequality (see \cite{Miclo} and \cite[Theorem 4.1]{W3}),
so does $\E_Q$ if and only if (\ref{rw5.3}) holds (see also
\cite{Chen}).

{\bf (iv)} follows from Corollary \ref{rwC4.3} and the facts that
$\E_0^{(1)}$ satisfies the super log-Sobolev inequality if and
only if $\theta>2$ (see \cite[Corollaris 2.5 and 3.3]{W3}), and so
does $\E_Q$ if and only if (\ref{rw1.7}) holds (see \cite{Mao2}).

 Finally,  by \cite[Corollary 2.5]{W3}, (\ref{rwSP0}) holds
 with $\bb_0(r)= \exp[c_0(1+r^{-1/\aa})]$ for some $c_0>0$
 if and only if $\theta\ge 2/(2-\aa)$.
Then by \cite[Corollary 1.2]{W4}, $\ll_\phi(\E_0^{(1)})>0$ for the
above $\phi$ if and only if   $\theta\ge 2/(2-\aa).$ Therefore,
{\bf (iii)} follows from Theorem \ref{rwT4.1}, Theorem
\ref{rwT4.2} (2) and and the fact that $\ll_\phi(\E_Q)>0$
if and only if (\ref{rw1.9}) holds. The
proof of this fact is  similar to that presented in
\cite{Mao1} for the log-Sobolev inequality, the only difference is
to use the so-called $N$-function
$\Psi(r):=|r|\{\log(1+|r|)\}^\aa$ in place of $|r|\log (1+|r|)$,
see \cite{FW} for details.
\end{proof}

\paragraph{Example 5.2. (with interacting diffusions)}
Let $M$ be compact and $\scr J:=\{J_A: A\subset\subset \N\}$ a
smooth potential with finite range; that is,
$J_A\in C^\infty(M^A)$ and vanishes if the diameter of $A$ is big enough.
A probability measure $\mu$ on $M^\Z$ is
called a Gibbs state with potential $\scr J$ if for any
$A\subset\subset \N$, its regular conditional distribution given
$x_{A^c}\in M^{A^c}$ is

$$\mu_{A|x}(\d y_A) := \ff 1 {Z_A(x_{A^c})} \exp[-U_A(y_A\times
x_{A^c})] \ll^A(\d y_A),$$ where $\ll^A$ is the volume measure on
$M^A$, $U_A:=\sum_{\LL: \LL\cap A\ne\emptyset} J_\LL$ and
$Z_A(x_{A^c})$ is the nomalization.
Let $\mu^{(n)}$ be the projection of $\mu$ on $M^{\{1,\cdots,n\}}$, and let
$\E^{(n)}$ be determined by (\ref{rwN}) with

$$\E_0(f,g):=\int_{M^\N} \sum_{k=1}^\infty\<\nn_k f,\nn_k g\>
\d\mu,\ \ \ f,g\in \scr FC^1(M^\N):= \bigcup_{A\subset\subset \N}
C^1(M^A),$$ where $\nn_k$ is the gradient w.r.t. the $k$-th
component. Assume that $\E_0$ satisfies the  log-Sobolev inequality

\beq\label{LLS}\mu(f^2\log f^2)\le c \E_0(f,f),\ \ \ f\in \scr FC^1(M^\N),
\mu(f^2)=1\end{equation} for some constant $c>0$. See e.g. \cite[Theorem
2.17]{DeS} for an explicit condition on $\scr J$ for (\ref{LLS})
to hold. Thus, if moreover
(\ref{rw1.9}) holds for $\aa=1$ so that $\E_Q$ satisfies the
log-Sobolev inequality, then Theorem \ref{rwT4.2} implies the log-Sobolev  for
$\E^{\GG_0}$.

\paragraph{Example 5.3. (the pure reaction case with finite $E$)}
Let $E=\{1,2,\cdots, N\}$ for some $N\ge 2$ and let $\mu$ be a
product probability measure on $E^\N$. Assume that
$\sum_{n=1}^\infty \vrr_n \e^{\dd n}<\infty$ for some $\dd>0$.
Then $\E_R^{\GG_0}$ satisfies the log-Sobolev inequality (i.e.
$\mathbf L(\E_R^{\GG_0})>0)$ if and only if (\ref{rw1.9}) holds with
$\aa=1$, while (\ref{rwSLS}) holds for $\E_R^{\GG_0}$ in place of
$\E^{\GG_0}$ for some $\bb$ if and only if (\ref{rw1.7}) holds.
Indeed, if $E$ is finite then the trivial Dirichlet form
$\E_0^{(1)}:=0$ satisfies (\ref{rwSLS0}) for some $\bb_0$, so that
the above assertions follow from Theorem \ref{rwT4.2} and
Corollary \ref{rwC4.3}.

\paragraph{Acknowledgement.} The authors would like to thank
Professor Mu-Fa Chen for useful suggestions and the referee for
careful comments.

\beg{thebibliography}{99}

\bibitem{AKR1} S. Albeverio, Yu. G. Kondratiev and M. R\"ockner, \emph{Analysis and geometry on configuration spaces,} J. Funct. Anal. 154(1998), 444--500

\bibitem{AKR2} S. Albeverio, Yu. G. Kondratiev and M. R\"ockner, \emph{Analysis and geometry on configuration spaces: the Gibbsian case,} J. Funct. Anal. 157(1998), 242--291.

\bibitem{BCC} L. Bertini, N. Cancrini and F. Cesi, \emph{The spectral gap for a Glauber-type
dynamics in a continuous gas,} Ann. Inst. H. Poincar\'e Probab. Statist. 38(2002), 91--108.

\bibitem{BL} S. Bobkov and M. Ledoux, \emph{On modified logarithmic Sobolev inequalities for Bernoulli
and Poisson measures,} J. Funct. Anal.
156(1998), 347--365.


\bibitem{CGT} J. Cheeger, M. Gromov and M. Taylor, \emph{Finite propagation speed, Laplace operator,
and geometry of complete Riemannian manifolds,} J. Diff. Geom. 17(1982), 15--53.

\bibitem{Chenbook} M.-F. Chen, \emph{From Markov Chains to
Non-Equilibrium Particle Systems}, World Scientific, Singapore,
1992.

\bibitem{Chen} M.-F. Chen,  \emph{Ergodic convergence rates of Markov processes -- eigenvalues,
inequalities and ergodic theory,} Proceedings of ICM (Beijing 2002), Vol. III, 41--52, Chinese High Edu.
Press,
Beijing 2002.

\bibitem{CW} M.-F. Chen and F.-Y. Wang, \emph{Cheeger's inequalities for general symmetric forms
and existence criteria for spectral gap,} Ann. Probab. 28(2000), 235--257.



\bibitem{DS} E. B. Davies and B. Simon, \emph{Ultracontractivity and the heat kernel
for Schr$\ddot{o}$dinger operators and Dirichlet Laplacians,} J.
Funct. Anal. 59(1984), 335-395.


\bibitem{DeS} J.-D. Deuschel and D. W. Stroock, \emph{Hypercontractivity and spectral gap
of symmetric diffusions with applications to the stochastic Ising
models,} J. Funct. Anal. 92(1990), 30--48.

\bibitem{GW} F.-Z. Gong and F.-Y. Wang, \emph{Functional inequalities for uniformly integrable
semigroups and application to essential spectrum,} Forum Math. 14(2002), 293--313.

\bibitem{Gross} L. Gross, \emph{Logarithmic Sobolev inequalities,}
Amer. J. Math. 97(1976), 1061--1083.

\bibitem{HS} R. A. Holley and D. W. Stroock, \emph{Nearest neighbor birth and death
processes on the real line,} Acta Math. 140(1987), 103--154.


\bibitem{Kond} Y. Kondratiev and E. Lytvynov, \emph{Glauber dynamics of
continuous particle systems,}  to appear in Ann. Inst. H. Poincar\'e Probab. Statist.

\bibitem{LO} R. Lata{\l}a and K. Oleszkiewicz, \emph{Between Sobolev and
Poincar\'e,} Lecture Notes Math. 1709, pp. 120-216, 1999.

\bibitem{Ledoux} M. Ledoux, \emph{On Talagrand's deviation inequalities for product measures,
}ESAIM: Probability and Statistics, 1(1996), 63-87.

\bibitem{MR} Z.-M. Ma and M. R\"ockner, \emph{Introduction to the
Theory of (Non-Symmetric) Dirichlet Forms,} Springer-Verlag, 1992.

\bibitem{Mao1} Y.-H. Mao, \emph{Logarithmic Sobolev inequalities for birth-death
 process and diffusion process on
the line,} Chinese J. Appl. Probab. Statist. 18(2002), 94-100.

\bibitem{Mao2} Y.-H. Mao, \emph{On supercontractivity for Markov semigroup,}
preprint.

\bibitem{Miclo} L. Miclo, \emph{An example of application of discrete
Hardy's inequalities,} Markov proc. Relat. Fields, 5(1999), 319--330.


\bibitem{R} M. R\"ockner, \emph{Stochastic analysis on configuration spaces: basic ideas and
 recent results,} In ``New Directions in Dirichlet Forms", AMS/IP
 Stud. Adv. Math. Vol. 8, pp. 157-231,  Amer. Math. Soc. Providence,
RI, 1998.

\bibitem{Preston} C. Preston, \emph{Spatial birth-and-death processes,} Proceedings of the 40th Session of
the International Statistical Institute (Warsaw 1975), Vol 2; Bull. Inst. Internat. Stat. 46(1975),
371--391.

\bibitem{RW1} M. R\"ockner and F.-Y. Wang, \emph{Weak
Poincar\'e inequalities and $L^2$-convergence rates of Markov
semigroups}, J. Funct. Anal. 185(2001), 564-603.

\bibitem{RW2} M. R\"ockner and F.-Y. Wang, \emph{On the spectrum of a class of (nonsymmetric) diffusion
operators,}   Bull. Lond. Math. Soc. 36(2004), 95--104.

 \bibitem{SZ} D. W. Stroock and B. Zegarlinski, \emph{The
 equivalence of the logarithmic Sobolev inequality and the
 Dobrushin-Shlosman mixing condition,} Comm. Math. Phys.
 144(1992), 303--323.

\bibitem{FW} F. Wang, \emph{ Lata{\l}a-Oleszkiewicz inequalities,} Ph.D. Thesis, Departmnent of Mathematics,
Beijing Normal University, 2003.

\bibitem{W11} F.-Y. Wang, \emph{Existence of spectral gap for elliptic operators,} Arkiv
Mat. 37(1999), 395--407.



\bibitem{W3} F.-Y. Wang, \emph{Functional inequalities for empty essential
spectrum,} J. Funct. Anal. 170(2000), 219-245.

\bibitem{W00b} F.-Y. Wang, \emph{Functional inequalities, semigroup properties and
spectrum estimates,} Infin. Dimens. Anal. Quant. Probab. Relat. Topics, 3(2000), 263--295.

\bibitem{W12} F.-Y. Wang, \emph{Logarithmic Sobolev inequalities: conditions and counterexamples,}
J. Operator Theory 46(2001), 183--197.

\bibitem{W10} F.-Y. Wang, \emph{Coupling, convergence rates of Markov processes and
weak Poincar\'e inequalities,} Sci. Sin. (A) 45(2002), 975--983.

\bibitem{W4} F.-Y. Wang, \emph{A generalization of Poincar\'e and log-Sobolev
inequalities},  Potential Analysis 22(2005), 1--15.

\bibitem{Wu1} L. Wu. \emph{A new modified logarithmic Sobolev inequality for
Point processes and several applications,} Probab. Theory Relat.
Fields 118(2000), 427-438.

\bibitem{Wu2} L. Wu, \emph{Estimates of spectral gap for continuous gas,}
 Ann. Inst. H. Poincar\'e Probab. Statist. 40(2004), 387--409.

\end{thebibliography}
\end{document}